\documentclass[12pt]{article}
\usepackage{amsmath,amsthm}
\usepackage{amsfonts,amsmath}
\usepackage{amssymb}
\usepackage{fullpage}
\usepackage{epsf}
\usepackage{color}
\usepackage{graphicx}
\usepackage{pdflscape}
\usepackage{rotating}

\usepackage{mathrsfs}
\usepackage{pgf,tikz}
\usepackage{multirow}

\newtheorem{theorem}{Theorem}[section]
\newtheorem{lemma}[theorem]{Lemma}

\newtheorem{corollary}[theorem]{Corollary}

\newtheorem{question}[theorem]{Question}
\theoremstyle{definition}
\newtheorem{definition}[theorem]{Definition}

\renewcommand{\geq}{\geqslant}
\renewcommand{\leq}{\leqslant}

\def\Z{\mathbb{Z}}

\def \md#1{{\,({\rm mod}\ #1)}}

\def \cF {{\cal F}}

\def \Z {\mathbb Z}

\def \R {\mathbb R}

\def\Z{\mathbb{Z}}

\definecolor{darkgreen}{RGB}{34, 150, 65}

\newcommand{\comment}[1]{}

\title{
Mutually orthogonal cycle systems
}

\author{Andrea C.~Burgess\thanks{Department of Mathematics and Statistics, University of New Brunswick, Saint John, NB, E2L~4L5, Canada.
\texttt{andrea.burgess@unb.ca}}
\and
Nicholas J.~Cavenagh\thanks{Department of Mathematics, The University of Waikato, Private Bag 3105, Hamilton 3240, New Zealand.
\texttt{nickc@waikato.ac.nz}}
\and
David A.~Pike\thanks{Department of Mathematics and Statistics, Memorial University of Newfoundland, St.~John's, NL, A1C~5S7, Canada.
\texttt{dapike@mun.ca}}
}

\begin{document}

\maketitle

\begin{abstract}
An ${\ell}$-cycle system ${\mathcal F}$ of a graph $\Gamma$ is a set of ${\ell}$-cycles which partition the edge set of $\Gamma$.  Two such cycle systems ${\mathcal F}$ and ${\mathcal F}'$ are said to be {\em orthogonal} if no two distinct cycles from ${\mathcal F}\cup {\mathcal F}'$ share more than one edge.  Orthogonal cycle systems naturally arise from face $2$-colourable polyehdra and in higher genus from Heffter arrays with certain orderings.  A set of pairwise orthogonal $\ell$-cycle systems  of $\Gamma$ is said to be a set of mutually orthogonal cycle systems of $\Gamma$.

Let $\mu(\ell,n)$ (respectively, $\mu'(\ell,n)$) be the maximum integer $\mu$ such that there exists a set of $\mu$ mutually orthogonal (cyclic) $\ell$-cycle systems of the complete graph $K_n$. We show that if $\ell\geq 4$ is even and $n\equiv 1\pmod{2\ell}$, then $\mu'(\ell,n)$, and hence $\mu(\ell,n)$, is bounded below by a constant multiple of $n/\ell^2$. In contrast, we obtain the following upper bounds: $\mu(\ell,n)\leq n-2$;
$\mu(\ell,n)\leq (n-2)(n-3)/(2(\ell-3))$ when $\ell \geq 4$;
$\mu(\ell,n)\leq 1$ when $\ell>n/\sqrt{2}$;
 and
 $\mu'(\ell,n)\leq n-3$ when $n \geq 4$.
We also obtain computational results for small values of $n$ and $\ell$.
\end{abstract}

\textbf{Keywords and MSC Code:} Orthogonal cycle decompositions; cyclic cycle systems; Heffter arrays; completely-reducible; super-simple. 05B30

\section{Introduction}
\label{Section:Intro}

We say that a graph $\Gamma$ {\em decomposes} into subgraphs $\Gamma_1, \Gamma_2, \ldots, \Gamma_t$, if the edge sets of the $\Gamma_i$ partition the edges of $\Gamma$. If
$\cF=\{\Gamma_i\mid 1\leq i\leq t\}$ where  $\Gamma_i \cong H$ for each $1\leq i\leq t$, then we say that $\cF$ is an $H$-decomposition of $\Gamma$.
An {\em $\ell$-cycle system} of a graph $\Gamma$ is a decomposition of $\Gamma$
into $\ell$-cycles. In the case where $\Gamma$ is the complete graph $K_n$ we say that there is an $\ell$-cycle system of {\em order} $n$.
Necessary and sufficient conditions for the existence of an $\ell$-cycle system of order $n$ were given in~\cite{AlspachGavlas2001,Sajna2002}; see also~\cite{BurattiRotational}.  Namely, at least one $\ell$-cycle system of order $n > 1$ exists if and only if $3 \leq \ell \leq n$, $n(n-1) \equiv 0 \pmod{2\ell}$ and $n$ is odd. 

Two $\ell$-cycle systems $\cF$ and $\cF'$ of the same graph $\Gamma$ are said to be {\em orthogonal} if, for all cycles $C\in \cF$ and $C'\in \cF'$, $C$ and $C'$ share at most one edge.
A set of pairwise orthogonal $\ell$-cycle systems
of $\Gamma$ is said to be a set of {\em mutually orthogonal} cycle systems of $\Gamma$.
In this paper we are interested in the maximum $\mu$ such that there exists a set of $\mu$ mutually orthogonal $\ell$-cycle systems of order $n$;
we denote this value by $\mu(\ell,n)$.


In the array below we exhibit a set of four mutually orthogonal cycle systems of order $9$. We have determined computationally that $\mu(4,9)=4$; i.e., this set is maximum.
{\small
$$\begin{array}{l}
  \{(1, 2, 3, 4),  (1, 3, 6, 5),  (1, 6, 2, 7),  (1, 8, 2, 9),  (2, 4, 7, 5),  (3, 5, 8, 7),  (3, 8, 6, 9),  (4, 5, 9, 8),  (4, 6, 7, 9)\},
\\[0.2ex]

  \{(1, 2, 6, 8),  (1, 3, 5, 7),  (1, 4, 8, 5),  (1, 6, 5, 9),  (2, 3, 6, 4),  (2, 5, 4, 9),  (2, 7, 3, 8),  (3, 4, 7, 9),  (6, 7, 8, 9)\},
\\[0.2ex]

  \{(1, 2, 8, 7),  (1, 3, 4, 6),  (1, 4, 9, 8),  (1, 5, 3, 9),  (2, 3, 8, 5),  (2, 4, 5, 6),  (2, 7, 5, 9),  (3, 6, 9, 7),  (4, 7, 6, 8)\},
\\[0.2ex]

  \{(1, 2, 9, 3),  (1, 4, 6, 9),  (1, 5, 4, 7),  (1, 6, 3, 8),  (2, 3, 7, 6),  (2, 4, 8, 7),  (2, 5, 6, 8),  (3, 4, 9, 5),  (5, 7, 9, 8)\}.
\end{array}
$$
}

Orthogonal cycle systems arise from face $2$-colourable embeddings of graphs on surfaces, which satisfy two conditions natural to polyhedra and similar phenomena: each pair of faces share at most one edge and each edge belongs to exactly two faces.

Let $\mu K_n$ be the multigraph in which each edge of $K_n$ is replaced by $\mu$ parallel edges.
A decomposition ${\mathcal F}$ of $\mu K_n$ into a subgraph $H$ is said to be {\em super-simple} if no two copies of $H$ share more than one edge, and {\em completely-reducible} if
${\mathcal F}$ partitions into $\mu$ decompositions of $K_n$. It follows that a set of $\mu$ mutually orthogonal cycle systems of $K_n$ is equivalent to a completely-reducible super-simple decomposition of
$\mu K_n$ into cycles; see \cite{CY2} for more details.

In the case $\ell=3$, observe that a pair of $\ell$-cycle systems is orthogonal if and only if the cycle systems are disjoint.
It is not hard to see that there are at most $n-2$ pairwise disjoint triple systems of order $n$; a set of systems which meets this bound is called a {\em large set} of disjoint Steiner triple systems, or LTS$(n)$.
An LTS$(7)$ does not exist \cite{Cayley}; however in \cite{LuJX1,LuJX2}, it is shown that an LTS$(n)$ exists if and only if $n>7$ and $n\equiv 1$ or $3\pmod{6}$, except for a finite list of possible exceptions. The exceptional cases are all solved in~\cite{Te}.

In this paper, we are often interested in {\em cyclic} cycle systems of the complete graph $K_n$.  Let $G$ be an additive group of order $n$ and suppose $K_n$ has vertex set $G$.  Given a cycle $C=(c_0,c_1,\ldots,c_{\ell-1})$ in $K_n$, for each element $g \in G$, define the cycle $C+g=(c_0+g, c_1+g, \ldots, c_{\ell-1}+g)$.  We say that a cycle system $\cF$ of $K_n$ is {\em $G$-regular} if, for any $C \in \cF$ and $g \in G$, we have that $C+g \in \cF$.  In the case that $G$ is a cyclic group, we refer to a $\Z_n$-regular cycle system as {\em cyclic}.  In a cyclic cycle system $\cF$, the {\em orbit} of the cycle $C \in \cF$ is the set of cycles $\{C+g \mid g \in \Z_n\}$; a cyclic cycle system can be completely specified by listing a set of {\em starter cycles}, that is, a set of representatives for the orbits of the cycles under the action of $\Z_n$.


The existence problem for cyclic cycle systems has attracted much attention.  Clearly, in order for a cyclic $\ell$-cycle system of odd order $n$ to exist, we must have that $3 \leq \ell \leq n$ and $\ell$ divides $n(n-1)/2$.  However, additional conditions for existence also come into play.  There is no cyclic $\ell$-cycle system of order $n$ when $(\ell,n)\in \{(3,9), (15,15)\}$; $\ell=n=p^m$ for some prime $p$ and integer $m \geq 2$; or $\ell < n < 2\ell$ and $\gcd(\ell,n)$ is a prime power~\cite{Buratti2004, BurattiDelFra2004}.  Buratti~\cite{Buratti2004} has conjectured that a cyclic $\ell$-cycle system of order $n$ exists for any other admissible pair $(\ell,n)$; this conjecture is still open.  The existence problem for cyclic cycle systems of the complete graph has been solved in a number of cases, including when $n \equiv 1$ or
$\ell\pmod{2\ell}$~\cite{BDF, BurattiDelFra2004, Kotzig, Rosa, Vietri} (see also~\cite{BlincoElZanatiVandenEynden, BryantGavlasLing, FuWu2004}), $\ell \leq 32$~\cite{WuFu}, $\ell$ is twice or thrice a prime power~\cite{Wu2,WuFu}, or $\ell$ is even and $m > 2\ell$~\cite{Wu}.


We explore the maximum $\mu'$ such that there exists a set of $\mu'$ mutually orthogonal {\em cyclic} $\ell$-cycle systems of order $n$;
this value is denoted by $\mu'(\ell,n)$.  Pairs of orthogonal cyclic cycle systems of the complete graph arise from Heffter arrays with certain orderings.
A {\em Heffter array $H(n;k)$} is an $n\times n$ matrix such that each row and column contains $k$ filled cells,
each row and column sum is divisible by $2nk+1$ and either $x$ or $-x$ appears in the array for each integer $1\leq x\leq nk$.
A Heffter array is said to have a {\em simple ordering} if, for each row and column, the entries may be cyclically ordered so that all partial sums are distinct modulo $2nk+1$.
The following was first shown by Archdeacon~\cite{A} as part of a more general result; consult~\cite{BCDY} to see this result stated more explicitly.

\begin{theorem}
If $H(n;k)$ is a Heffter array with a simple ordering, then there exists a pair of orthogonal cyclic decompositions of $K_{2nk+1}$ into $k$-cycles. In particular, $\mu'(k,2nk+1)\geq 2$.
\end{theorem}

Thus the following is implied by existing literature on Heffter arrays.
\begin{theorem} {\rm \cite{ADDY,BCDY,CMPP,DW}}
 Let $n\geq k$. Then
 $\mu'(k,2nk+1)\geq 2$ whenever:
\begin{itemize}
\item $k\in \{3,5,7,9\}$ and $nk\equiv 3\pmod{4}$; 
\item $k\equiv 0\pmod{4}$; 
\item $n\equiv 1\pmod{4}$ and $k\equiv 3\pmod{4}$; 
\item $n\equiv 0\pmod{4}$ and $k\equiv 3\pmod{4}$  (for large enough $n$). 
\end{itemize}
\end{theorem}

With an extra condition on the orderings of the entries of a Heffter array, these orthogonal cycle systems in turn biembed to yield a face $2$-colourable embedding on an orientable surface. Face $2$-colourable embeddings on orientable surfaces have been studied for a variety of combinatorial structures \cite{DM,GG,GrM,GM}.  Recently, Costa, Morini, Pasotti and Pellegrini~\cite{CMPP2020} employed a generalization of Heffter arrays to construct pairs of orthogonal $\ell$-cycle systems of the complete multipartite graph in certain cases.

In \cite{CY2}, it is shown that for every graph $H$ and fixed integer $k\geq 1$, for sufficiently large $n$ (satisfying some elementary necessary divisibility conditions), there exists a set of $k$ pairwise orthogonal decompositions of $K_n$ into $H$ (i.e., no two copies of $H$ share more than one edge).
Aside from this quite general asymptotic result, to our knowledge, sets of mutually orthogonal $\ell$-cycle systems of size greater than $2$ have not been studied for $\ell\geq 4$.

In this paper, our focus for cyclic cycle systems is in the case $n \equiv 1\pmod{2\ell}$, for which it is possible to construct a cyclic $\ell$-cycle system with no short orbit.
In particular, we will find lower bounds on $\mu(\ell,n)$ by constructing sets of mutually orthogonal cyclic even cycle systems.
Specifically, we show that if $\ell$ is even and $n \equiv 1\pmod{2\ell}$, then $\mu'(\ell,n)$ is bounded below by a constant multiple of $n/\ell^2$, i.e., $\mu'(\ell,n) = \Omega(n/\ell^2)$.
 Our main result is as follows.
\begin{theorem}\label{MainTheorem}
If $\ell \geq 4$ is even, $n \equiv 1\pmod{2\ell}$ and $N=(n-1)/(2\ell)$, then
\[
\mu(\ell,n) \geq \mu'(\ell,n) \geq \frac{N}{a\ell+b}-1,
\]
where
\[
(a,b) = \left\{
\begin{array}{ll}
(4,-2), & \mbox{if } \ell \equiv 0\pmod{4}, \\
(24,-18), & \mbox{if } \ell \equiv 2\pmod{4}.
\end{array}
\right.
\]
\end{theorem}

In Section~\ref{Section:MO4CS}, when $\ell=4$, we improve the bound of Theorem~\ref{MainTheorem} to $\mu(\ell,n) \geq \mu'(\ell,n) \geq 4N$ (Lemma~\ref{CycleLength4}).
Section~\ref{Section:Prelim} establishes some notation and preliminary results.
The general result for $\ell \equiv 0 \pmod{4}$ is proved in Section~\ref{Section:4k} (Theorem~\ref{case8k}), while the bound for $\ell \equiv 2\pmod{4}$ is proved in Section~\ref{Section:4k+2} (Theorem~\ref{Case4k+2}).
In contrast, in Section~\ref{conclus} we establish upper bounds, namely
$\mu(\ell,n)\leq n-2$;
$\mu(\ell,n)\leq (n-2)(n-3)/(2(\ell-3))$ for $\ell\geq 4$;
$\mu(\ell,n)\leq 1$ for $\ell>\sqrt{n(n-1)/2}$;
 and $\mu'(\ell,n)\leq n-3$ for $n\geq 4$.
Finally, computational results for small values are given in the appendix.

\section{Mutually orthogonal 4-cycle systems}
\label{Section:MO4CS}

Clearly $n\equiv 1\md{8}$ is a necessary condition for a  decomposition  of $K_n$ into $4$-cycles, cyclic or otherwise.
Let $[a,b,c,d]_n$ denote the $\Z_n$-orbit of the $4$-cycle $(0,a,a+b,a+b+c)$, where $a+b+c+d$ is divisible by $n$.
Observe that $[a,b,c,d]_n=[-d,-c,-b,-a]_n$.
Where the context is clear,  we write $[a,b,c,d]_n=[a,b,c,d]$.
Let $D_n=\{1,2,\dots , (n-1)/2\}$; that is, $D_n$ is the set of {\em differences} in $\Z_n$.
We consider $\Z_n$ as the set $\pm D_n \cup \{0\}$.

By observation, the maximum size of a set of mutually orthogonal cyclic $4$-cycle systems of $K_9$ is $\mu'(4,9)=2$.
Two such systems are  $[1,-2,4,-3]_9$ and $[1,-3,4,-2]_9$.
In the non-cyclic case, an exhaustive computational search indicates that the maximum size of a set of mutually orthogonal $4$-cycle systems of $K_9$ is
$\mu(4,9)=4$; see the example given in Section~\ref{Section:Intro}.

\begin{lemma}\label{CycleLength4}
If $n\equiv 1\md{8}$ and $n\geq 17$, then there exists a set of $(n-1)/2$ mutually orthogonal cyclic $4$-cycle systems of order $n$.
In particular, $\mu'(4,n)\geq (n-1)/2$.
\end{lemma}

\begin{proof}
We first describe how to construct a set of $(n-5)/2$ mutually orthogonal cyclic $4$-cycle systems; then we add two more by making some adjustments.

Let $N=(n-1)/8$.
For each $i,j$ with $1\leq i<j\leq 2N$,
let $C_{i,j}$ and $C_{i,j}'$ be the pair of orbits of $4$-cycles:
 $$C_{i,j}:=\{[2i-1,2j,-2i,-(2j-1)]\},\quad C_{i,j}':=\{[2i-1,-(2j-1),-2i,2j]\}.$$
Next, let $F_1,F_2,\dots F_{2N-1}$ be a set of $1$-factors which decompose the complete graph on vertex set $\{1,2,\ldots,2N\}$.

For each $1$-factor $F_k$,
the sets $$\cF_k:=\mathop{\bigcup_{\{i,j\}\in F_k}}_{i<j} C_{i,j}
\quad \mbox{ and } \quad
\cF_k':=\mathop{\bigcup_{\{i,j\}\in F_k}}_{i<j} C_{i,j}'
$$ each describe a cyclic decomposition of $K_n$ into $4$-cycles.
Observe that the set of such decompositions constitutes a mutually orthogonal set of size $4N-2=(n-5)/2$.

We next make an adjustment to extend this set.
Without loss of generality, let $F_1=\{\{1,2\},\{3,4\},\dots ,\{2N-1,2N\}\}$.
Replace $\cF_1$ and $\cF_1'$ with:

$$\begin{array}{l}
\cF_{\ast}=\{[4i-3,-(4i-2),-(4i-1),4i]\mid 1\leq i\leq N\}, \\
\cF_{\ast}'= \{[4i-3,4i,-(4i-1),-(4i-2)]
\mid 1\leq i\leq N\}.
\end{array}$$
Then, we can add another pair of cyclic decompositions, orthogonal to each
decomposition in $\{\cF_{\ast},\cF_{\ast}',\cF_2,\ldots,\cF_{2N-1},\cF_2',\ldots,\cF_{2N-1}'\}$, given by:
$$\cF_{2N}:=\{[1,-3,4N,-(4N-2)]\}\cup \{[4i+1,-(4i+3),4i,-(4i-2)]\mid 1\leq i <N\}$$
and
$$\cF_{2N}':=\{[1,-(4N-2),4N,-3]\}\cup \{[4i+1,-(4i-2),4i,-(4i+3)]\mid 1\leq i <N\}.$$
(Note that orthogonality requires $N\geq 2$ at this final step.)
\end{proof}

In the case $n=17$, we have computationally determined that $\mu'(4,17)=10$, which improves on the bound given in Lemma~\ref{CycleLength4}.  
A set of ten mutually orthogonal cyclic 4-cycle systems of order 17 is given in the appendix.

We exhibit the methods of the previous proof in the case $n=25$. We start with a $1$-factorization of $K_6$:
$$
\begin{array}{l}
F_1 = \{\{1,2\}, \{3,4\}, \{5,6\}\}, \\
F_2 = \{\{1,3\}, \{2,6\}, \{4,5\}\}, \\
F_3 = \{\{1,4\}, \{2,5\}, \{3,6\}\}, \\
F_4 = \{\{1,5\}, \{2,3\}, \{4,6\}\}, \\
F_5 = \{\{1,6\}, \{2,4\}, \{3,5\}\}.
\end{array}
$$
The resulting $12$ mutually orthogonal cyclic $4$-cycle systems of order $25$ are given by:

$$\begin{array}{l}
\cF_{\ast}=\{[1,-2,-3,4], [5,-6,-7,8], [9,-10,-11,12]\}, \\
\cF_{\ast}'=\{[1,4,-3,-2], [5,8,-7,-6], [9,12,-11,-10]\},\\

\cF_{2}= \{[1,6,-2,-5], [3,12,-4,-11], [7,10,-8,-9]\}, \\
\cF_{2}'= \{[1,-5,-2,6], [3,-11,-4,12], [7,-9,-8,10]\}, \\

\cF_{3}= \{[1,8,-2,-7], [3,10,-4,-9], [5,12,-6,-11]\}, \\
\cF_{3}'=\{[1,-7,-2,8], [3,-9,-4,10], [5,-11,-6,12]\}, \\

\cF_{4}= \{[1,10,-2,-9], [3,6,-4,-5], [7,12,-8,-11]\}, \\
\cF_{4}'=\{[1,-9,-2,10], [3,-5,-4,6], [7,-11,-8,12]\}, \\

\cF_{5}= \{[1,12,-2,-11], [3,8,-4,-7], [5,10,-6,-9]\}, \\
\cF_{5}'=\{[1,-11,-2,12], [3,-7,-4,8], [5,-9,-6,10]\}, \\

\cF_{6}=\{[1,-3,12,-10], [5,-7,4,-2], [9,-11,8,-6]\}, \\
\cF_{6}'= \{[1,-10,12,-3], [5,-2,4,-7], [9,-6,8,-11]\}.
\end{array}$$

Through computational means we determined that this collection of 12 mutually orthogonal cyclic 4-cycle systems of order 25 is maximal.
However, it is not maximum, as we also established computationally that $\mu'(4,25) \geq 17$.

\section{Preliminary lemmas for cycle length greater than 4}
\label{Section:Prelim}

In this section, we introduce notation and basic results which will be needed later to construct mutually orthogonal cycle systems with even cycle length $\ell \geq 6$.

Henceforth, for any integers $a$ and $b$ with $a\leq b$, $[a,b]$ is the set of integers $\{a,a+1,\dots ,b\}$. For $a,b\in \R$ with $a<b$, we also use the notation $(a,b)$ to denote the set of {\em integers} strictly between
$a$ and $b$.

Let the vertices of the complete graph $K_n$ be labelled with $[0,n-1]$, where $n$ is odd. Then the {\em difference} associated with an edge $\{a,b\}$ is defined to be the minimum value in the set
$\{|a-b\pmod {n}|,|b-a\pmod{n}|\}$. Let $e_1$ and $e_2$ be two edges of differences $d$ and $e$, respectively.
Then we may write $e_1=\{a,a+d\pmod{n} \}$ and $e_2=\{b,b+e\pmod{n} \}$, where $a,b\in [0,n-1]$ are uniquely determined.
We define the {\em distance} between
$e_1$ and $e_2$ to be the minimum value in the set $\{|a-b\pmod {n}|,|b-a\pmod{n}|\}$.
Given a cycle $C$ with vertices in $\Z_n$, the set $\Delta C$ is defined to be the multiset of differences of the edges of $C$.

The idea is to construct cyclic systems using so-called {\em balanced} sets of differences.  The following definitions and lemma appear in~\cite{BurgessMerolaTraetta}.

\begin{definition}
If $D=\{d_1, d_2, \ldots, d_{2k}\}$ is a set of positive integers, with $d_{i}< d_{i+1}$ for $i\in[1,2k-1]$,
the {\em alternating difference pattern} of $D$ is the sequence $(s_1, s_2, \ldots, s_k)$ where
$s_i=d_{2i}-d_{2i-1}$ for every $i\in[1,k]$. Furthermore, $D$ is said to be {\em balanced} if there exists an integer
$\tau\in[1,k]$ such that
$\sum_{i=1}^{\tau} s_i = \sum_{i=\tau+1}^k s_i$.
\end{definition}

\begin{definition}
Let $D=\{d_1, d_2, \ldots, d_{2k}\}$ be a balanced set of positive integers. Let $\delta_1$, $\delta_2, \dots , \delta_{2k}$ be the sequence obtained by reordering the integers in $D$ as
follows:
$$\delta_i =
\left\{
\begin{array}{ll}
d_i & \mbox{\rm if\ }1\leq i\leq 2\tau-1, \\
d_{i+1} & \mbox{\rm if\ }2\tau\leq i\leq 2k-1, \\
d_{2\tau} & \mbox{\rm if\ }i=2k.
\end{array}\right.$$
Set $c_0=0$ and $c_i=\sum_{h=1}^i (-1)^h\delta_h$ for $1\leq i\leq 2k-1$.
We then define $C(D):=(c_0,c_1,\dots ,c_{2k-1})$.
\end{definition}

\begin{lemma}\label{lemma1} {\rm (Lemma 3.2 of \cite{BurgessMerolaTraetta}).}
Let $k\geq 2$. If $D$ is a balanced set of $2k$ positive integers, then $C(D)$ is a
$2k$-cycle satisfying $\Delta C(D)=D$ and
vertex set $V(C(D))\subset [-d, d']$, where $d=\max D$ and $d'=\max (D\setminus\{d\})$.
\end{lemma}

\begin{corollary}\label{corollary1}
Let $k\geq 2$ and $n\equiv 1\pmod{4k}$. Suppose that the set
$[1,(n-1)/2]$ partitions into sets $D_1, D_2,\dots, D_{(n-1)/(2k)}$, each of which is balanced and of size $2k$.
Then cycles $C(D_i)$, $i \in [1, (n-1)/(2k)]$, form a set of starter cycles for a cyclic $2k$-cycle decomposition of $K_n$; in particular, the set
$$\{C(D_i)+j\mid i\in [1, (n-1)/2k], j\in [0,n-1]\}$$
is a cyclic decomposition of $K_n$ into $2k$-cycles.
\end{corollary}

\begin{proof}
Let $i\in [1, (n-1)/2k]$. Since $D_i\subset [1,(n-1)/2]$, Lemma \ref{lemma1} implies that  $V(C(D_i))\subset [-(n-1)/2, (n-1)/2]$. Thus the vertices of $V(C(D_i))$ are distinct in $\Z_n$. The result follows.
\end{proof}

Our general strategy will be to show that a pair of cyclically generated cycle systems is orthogonal by showing that the sets of differences from any two cycles in different orbits share at most one element. To this end, the following lemma will be used in Sections~\ref{Section:4k} and~\ref{Section:4k+2}.

\begin{lemma}
Let $\delta, N>0$ and let $d$ and $d'$ be integers such that $d,d'\in (N/2-\delta N, N/2+\delta N)$.
Let $\alpha, \alpha'$ be integers such that $1\leq \alpha < \alpha' \leq (1-2\delta)/4\delta$.
Then  $\alpha d < \alpha' d'$.
\label{squeezy}
\end{lemma}

\begin{proof}
For each positive integer $s$, define
$$I_s = \{si\mid N/2-\delta N <i < N/2+\delta N; i\in {\mathbb R}\}.$$
Let $m=\lfloor \frac{1-2\delta}{4\delta}\rfloor$, and let $S=[1,m]$.
Observe that $\alpha, \alpha'\in S$.
Now,  $\delta \leq 1/(4m+2)$ implies that:
\[
\begin{array}{rrcl}
 & m(1+2\delta) & \leq & (m+1)(1- 2\delta)  \\
\Rightarrow & m(N/2 + \delta N) & \leq & (m+1)(N/2 -  \delta N). \\
\end{array}
\]
It follows that for each $s\in S$, every element of $I_s$ is strictly less than every element of $I_{s+1}$.
Since $\alpha d\in I_{\alpha}$ and  $\alpha' d'\in I_{\alpha'}$, it follows that
 $\alpha d< \alpha' d'$.
\end{proof}

The following variation of Lemma~\ref{squeezy} will be used in Section~\ref{Section:4k+2}.

\begin{corollary}
Let $\delta, N>0$ and let $d$ and $d'$ be integers such that $d,d'\in (N/3-\delta N, N/3+\delta N)$.
Let $\alpha, \alpha'$ be integers such that $1\leq \alpha < \alpha' \leq (1-3\delta)/6\delta$.
Then  $\alpha d < \alpha' d'$.
\label{squeezier}
\end{corollary}

\begin{proof}
If $m$ is a positive integer,
$m\leq (1-3\delta)/6\delta$ implies that
$$m(N/3 + \delta N) \leq (m+1)(N/3 -  \delta N).$$
The remaining argument is similar to the previous lemma.
\end{proof}

\section{Orthogonal sets of $4k$-cycle systems with $k \geq 2$}
\label{Section:4k}

Our aim in this section is to prove Theorem~\ref{case8k}. In particular, for each $k \geq 2$ and $n \equiv 1\pmod{8k}$, we will show that $\mu'(n,4k) = \Omega(n/k^2)$. That is, we construct a set of mutually orthogonal $4k$-cycle decompositions of $K_n$ of size at least $cn/k^2$ where $c$ is a constant.

Let $N$ and $k$ be positive integers and let $n=8kN+1$.
For each integer $d\in (N/2-N/(16k-2), N/2)$,  we construct a cyclic $4k$-cycle decomposition of $K_n$ which we will denote by ${\mathcal F}(d)$.

The first $d$ starter cycles in $\mathcal{F}(d)$ use the set of differences $[1,4kd]$.  For $i \in [1,d]$, let
\[
S_{d,i}=\{i,d+i,2d+i,\ldots,(4k-1)d+i\}.
\]
Observe that the set $S_{d,i}$ is balanced, with $\tau=k$, for each $i \in [1,d]$.

Henceforth in this section, let $e:=N-d$.  (In effect, $e$ is a function of $d$.)  Observe that $e\in (N/2, N/2 + N/(16k-2))$.
The remaining $e$ starter cycles in  ${\mathcal F}(d)$ use differences $[4kd+1, 4kN]$.  For $i \in [1, e]$, take
\[
T_{e,i} = \{4kd+i, 4kd+e+i, 4kd+2e+i, \ldots, 4kd+(4k-1)e+i\}.
\]
Observe that the set $T_{e,i}$ is balanced for each $i\in [1,e]$, where $\tau=k$.
Moreover, since $4kd+4ke=4kN$, we have that
\[
\left(\bigcup_{i=1}^{d} S_{d,i}\right) \cup \left( \bigcup_{i=1}^{e} T_{e,i} \right)=[1,4kN],
\]
so by Lemma~\ref{corollary1},
the set of cycles
$${\mathcal F}(d):=\{C(S_{d,i})\mid i\in [1,d]\}\cup \{C(T_{e,i})\mid i\in [1,e]\}$$
is a set of starter cycles for a cyclic $4k$-cycle system of order $n=8kN+1$.

In order to show that we have constructed an orthogonal set of decompositions, we will make use of the following, which is a direct consequence of Lemma~\ref{squeezy}.

\begin{lemma}
Let
$d,d'\in (N/2-N/(16k-2), N/2)$ where $d\neq d'$ and let $e=N-d$ and $e'=N-d'$. Let
$\alpha, \alpha' \in [1,4k-1]$.
Then no two of $\alpha d$, $\alpha' d'$, $\alpha e$ and $\alpha' e'$ are equal.
 Moreover, if $\alpha < \alpha'$ then $\alpha d < \alpha' d'$ and $\alpha e < \alpha'e'$.
\label{Squeezy2}
\end{lemma}

\begin{lemma}
Let $d,d'\in (N/2-N/(16k-2), N/2)$ where $d\neq d'$.
Then
the decompositions ${\mathcal F}(d)$ and  ${\mathcal F}(d')$, as defined above, are orthogonal.
\end{lemma}

\begin{proof}
 In what follows,
  $d\neq d'$, $e=N-d$ and $e'=N-d'$. Observe that
   $e,e'\in (N/2, N/2 + N/(16k-2))$.

It suffices to show that if  $C$ is a cycle from $\mathcal{F}(d)$ and $C'$ is a cycle from  $\mathcal{F}(d')$, then
$C$ and $C'$ share at most one difference.
Equivalently, we will show that:
\begin{enumerate}
\item [{\bf (i)}:] {\em For any $i\in [1,d]$ and $i'\in [1,d']$, $|S_{d,i}\cap S_{d',i'}|\leq 1$};
\item  [{\bf (ii)}:] {\em For any $i\in [1,e]$ and $i'\in [1,e']$, $|T_{e,i}\cap T_{e',i'}|\leq 1$}; and
\item  [{\bf (iii)}:] {\em For any $i\in [1,d]$ and $i'\in [1,e']$, $|S_{d,i}\cap T_{e',i'}|\leq 1$}.
\end{enumerate}

To show (i), suppose to the contrary that $\{x,y\} \subseteq S_{d,i} \cap S_{d',i'}$ with $x<y$. Thus $y-x = \alpha d = \alpha' d'$ for some $\alpha, \alpha'
\in [1,4k-1]$, contradicting Lemma~\ref{Squeezy2}.  The justification of (ii) is similar.  For (iii), if $x,y \in S_{d,i}\cap T_{e',i'}$ with $x<y$, then $y-x = \alpha d$ for some $\alpha \in [1,4k-1]$ (since $x,y \in S_{d,i}$) and $y-x=\alpha' e'$ for some
$\alpha' \in [1,4k-1]$ (since $x,y \in T_{e',i'}$), so $\alpha d = \alpha'e'$, which again contradicts Lemma~\ref{Squeezy2}.
\end{proof}

Since $n=8Nk+1$, we have the following theorem.
\begin{theorem}
Let $k \geq 2$ and $n=8Nk+1$.  There is a set of mutually orthogonal cyclic $4k$-cycle systems of order $n$ of size at least
$$\frac{N}{16k-2}-1=\frac{n-1}{8k(16k-2)}-1.$$
Thus, if $n\equiv 1\pmod{8k}$,
$$\mu(n,4k)\geq \mu'(n,4k)\geq \frac{n-1}{8k(16k-2)}-1.$$
\label{case8k}
\end{theorem}

\section{Orthogonal sets of $(4k+2)$-cycles}
\label{Section:4k+2}

Let
$N$ and $k$ be positive integers and let
$n=2(4k+2)N+1$.  For each $d \equiv N\pmod{2}$ with $d\in (N/3 - N/(48k+15), N/3)$, we form a cyclic $(4k+2)$-cycle decomposition $\mathcal{F}(d)$ of $K_n$.  Let $e=(N-d)/2$, and observe that $N/3 < e < N/3 + N/(2(48k+15))<  N/3 + N/(48k+15)$.
Thus $e\in  (N/3,  N/3 + N/(48k+15))$.

For $i \in [1,d]$, let
\[
S_{d,i,1} = \{i, d+i, 2d+i, \ldots, (4k-1)d+i\} \mbox{ and } S_{d,i,2} = \{4kN+4e+i, (4k+2)N-i+1\},
\]
and let $S_{d,i} = S_{d,i,1} \cup S_{d,i,2}$.

Now, when constructing the cycles containing differences in $S_{d,i}$, instead of $(4k+2)N-i+1$, we will use the {\em negative} of this difference modulo $n$, that is, the value
\[
(8k+4)N+1 - ((4k+2)N-i+1) = (4k+2)N+i.
\]
We construct a starter cycle $C'(S_{d,i})$ using the set of differences $S_{d,i}$ but in a slightly different way to Lemma~\ref{lemma1}.
\begin{eqnarray*}
C'(S_{d,i}) &=& (0,-i,d,-d-i,\ldots, kd, -kd-i, \\
&& (k+2)d, -(k+1)d-i, (k+3)d, -(k+2)d-i, \ldots, 2kd, -(2k-1)d-i, \\
&& (4k+2)N-(2k+1)d, -(2k+1)d-i).
\end{eqnarray*}
(Note that in the case $k=1$, $C'(S_{d,i}) = (0,-i,d,-d-i,4N+e-d, -3d-i)$.)

\begin{lemma}
Let $i \in [1,d]$. Working modulo $n$, the ordered sequence $C'(S_{d,i})$ is a
 $(4k+2)$-cycle with difference set $S_{d,i}$.
\label{specialcycle}
\end{lemma}

\begin{proof}
To see that no vertices are repeated (modulo $n$) within the sequence $C'(S_{d,i})$, it suffices to observe that:
\[
\begin{array}{l}
-(4k+2)N < -(2k+1)d-i < -(2k-1)d-i < -(2k-2)d-i < \cdots < -d-i < -i \\
< 0 < d < 2d < \cdots < kd < (k+2)d < (k+3)d < \cdots < 2kd \\
< (4k+2)N-(2k+1)d < (4k+2)N.
\end{array}
\]
By inspection, and since $(4k+2)N-(2k+1)d = 4kN+4e-(2k-1)d$ and $n-((4k+2)N-i+1)=(4k+2)N+i$, the set of differences of the edges of the cycle $C'(S_{d,i})$ is $S_{d,i}$.
\end{proof}

Note that
\[
\bigcup_{i=1}^d S_{d,i} = [1,4kd] \cup [4kN+4e+1,4kN+4e+d] \cup [(4k+2)N-d+1,(4k+2)N];
\]
since $4kN+4e+d = (4k+2)N-d$, we have that
\[
\bigcup_{i=1}^d S_{d,i} = [1,4kd] \cup [4kN+4e+1,(4k+2)N].
\]

For $j,\ell \in [1,e]$, let
\[
\begin{array}{l}
T_{e,j,1} = \{4kd+j, 4kd+e+j, \ldots, 4kd+(4k-1)e+j\}, \\
T_{e,j,2} = \{ 4kN+j, 4kN+2e+j\}, \\
U_{e,\ell,1} = \{4kd+4ke+\ell, 4kd+(4k+1)e+\ell, \ldots, 4kd+(8k-1)e+\ell\}, \\
U_{e,\ell,2} = \{4kN+e+\ell, 4kN+3e+\ell\},
\end{array}
\]
and set $T_{e,j} = T_{e,j,1} \cup T_{e,j,2}$ and $U_{e,\ell} = U_{e,\ell,1} \cup U_{e,\ell,2}$.

The sets $T_{e,j}$ and $U_{e,\ell}$ are each balanced with $\tau=k+1$.  We have that
\[
\left(\bigcup_{j=1}^e T_{e,j}\right) \cup \left(\bigcup_{\ell=1}^e U_{e,\ell}\right) = [4kd+1,4kd+8ke] \cup [4kN+1,4kN+4e] = [4kd+1,4kN+4e],
\]
since $4kd+8ke=4kN$.  Observe that for fixed $d$,
\[
\left(\bigcup_{i=1}^d S_{d,i}\right) \cup \left(\bigcup_{j=1}^e T_{e,j} \right) \cup \left(\bigcup_{\ell=1}^e U_{e,\ell}\right) = [1,(4k+2)N],
\]
and thus by Lemmas~\ref{corollary1} and~\ref{specialcycle}, the set of cycles
\[
\mathcal{F}(d)=\{C'(S_{d,i}) \mid i \in [1,d]\} \cup \{C(T_{e,j}) \mid j \in [1,e]\} \cup \{C(U_{e,\ell}) \mid \ell \in [1,e]\}
\]
is a set of starter cycles for a $(4k+2)$-cycle decomposition of $K_n$.

In order to show that the decompositions $\mathcal{F}(d)$, $d \in (N/3-N/(48k+15), N/3)$, are orthogonal, we will make use of the following lemma which is directly implied by Corollary~\ref{squeezier}.

\begin{lemma} \label{squeezy4}
Let $d\neq d'$, $e\neq e'$ and $$d,d',e,e'\in \left(\frac{N}{3} - \frac{N}{48k+15}, \frac{N}{3} + \frac{N}{48k+15}\right).$$  Let $\alpha, \alpha' \in [1,8k+2]$.  Then $\alpha d \neq \alpha' d'$ and $\alpha e \neq \alpha' e'$.
 Moreover, if $\alpha < \alpha'$, then $\alpha d < \alpha' d'$ and $\alpha e < \alpha' e'$.
\end{lemma}

\begin{lemma}
Suppose that $\beta d + i = \beta'd'+i'$, where $\beta, \beta' \in [0, 4k-1]$,
$i\in [1,d]$, $i' \in [1,d']$ and $d' < d$.  Then either $\beta'=\beta$ or $\beta'=\beta+1$.
\label{prelim}
\end{lemma}

\begin{proof}
From Lemma~\ref{squeezy4}, $(\beta+1)d < (\beta+2)d'$.  Now,
\[
(\beta-1)d'+i' \leq \beta d' \leq  \beta d < \beta d+ i
\]
and
\[
\beta d + i \leq (\beta+1)d < (\beta+2)d' < (\beta+2)d'+i';
\]
hence
\[
(\beta-1)d'+i' < \beta d + i < (\beta+2)d'+i'.
\]
\end{proof}

\begin{lemma}
Let $d\neq d'$ such that $d,d' \equiv N\pmod{2}$ and
$$d,d'\in \left(\frac{N}{3} - \frac{N}{48k+15}, \frac{N}{3} + \frac{N}{48k+15}\right).$$
  Let $e=(N-d)/2$ and $e'=(N-d')/2$.  Let $i \in [1,d]$, $i' \in [1,d']$, $j, \ell \in [1,e]$ and $j', \ell' \in [1,e']$.  Then for each $X \in \{S_{d,i}, T_{e,j}, U_{e,\ell}\}$ and each $Y \in \{S_{d',i'}, T_{e',j'}, U_{e',\ell'}\}$, $|X \cap Y| \leq 1$ \em{with the exception} $S_{d,i} \cap S_{d',i} = \{i,(4k+2)N+i\}$.
\label{Intersections2mod4}
\end{lemma}

\begin{proof}
Recall from the start of this section that $e,e'\in (N/3, N/3 + N/(48k+15))$. In what follows, we frequently apply Lemma~\ref{squeezy4} to $d,d',e$ and $e'$.
To prove the lemma, it suffices to show the following:
\begin{enumerate}
\item[{\bf (i):}] $S_{d,i} \cap S_{d',i} = \{i,(4k+2)N-i+1\}$ and if $i \neq i'$ then $|S_{d,i} \cap S_{d',i'}|\leq 1$;
\item[{\bf (ii):}] $|T_{e,j} \cap T_{e',j'}| \leq 1$, $|U_{e,\ell} \cap U_{e',\ell'}| \leq 1$ and $|T_{e,j} \cap U_{e',\ell'}| \leq 1$;
\item[{\bf (iii):}] $|S_{d,i} \cap T_{e',j'}| \leq 1$ and $|S_{d,i} \cap U_{e',\ell'}| \leq 1$.
\end{enumerate}
\bigskip

\noindent {\bf Proof of (i):}  In this case, we may assume without loss of generality that $d'<d$.  We note that
\begin{center}
\begin{tabular}{l}
$4kN+4e'+i' > 4kN > 4kd \geq (4k-1)d+i$ \ and  \\
$4kN+4e+i > 4kN > 4kd' \geq (4k-1)d'+i'$,
\end{tabular}
\end{center}
so $S_{d,i,1} \cap S_{d',i',2} = S_{d',i',1} \cap S_{d,i,2} = \emptyset$.

Now, supposing that $|S_{d,i,1} \cap S_{d',i',1}| \geq 2$, it follows that for some $x$, $(x,x+\alpha d) = (x, x+\alpha' d')$ where $\alpha, \alpha' \in [1,4k-1]$; thus $\alpha d = \alpha' d'$, in contradiction to Lemma~\ref{squeezy4}.  Next, supposing that $|S_{d,i,2} \cap S_{d',i',2}| \geq 2$, then either
\begin{center}
\begin{tabular}{ll}
(a) & $4kN+4e+i=4kN+4e'+i'$ and $(4k+2)N-i+1 = (4k+2)N-i'+1$, or \\
(b) & $4kN+4e+i = (4k+2)N-i'+1$ and $4kN+4e'+i' = (4k+2)N-i+1$.
\end{tabular}
\end{center}
In both cases, it is straightforward to check that $e=e'$, a contradiction.

Thus if $|S_{d,i} \cap S_{d',i'}| \geq 2$, it must be that $|S_{d,i,1} \cap S_{d',i',1}|=1$ and $|S_{d,i,2} \cap S_{d',i',2}|=1$.
If $i=i'$ then $\{i,(4k+2)N-i+1\}\subseteq S_{d,i} \cap S_{d',i'}$.
Moreover, recalling that $S_{d,i,1} \cap S_{d',i',2} = S_{d',i',1} \cap S_{d,i,2} = \emptyset$, it follows that $|S_{d,i} \cap S_{d',i'}|=2$.  Hence if $i=i'$, then $S_{d,i} \cap S_{d',i} = \{i,(4k+2)N-i+1\}$.
We now assume that $i \neq i'$.  From Lemma~\ref{prelim}, $|S_{d,i,1} \cap S_{d',i',1}|=1$ implies that either
\begin{center}
\begin{tabular}{ll}
(a) & $\beta d + i = \beta d' + i'$, or \\
(b) & $\beta d + i = (\beta+1) d' + i'$
\end{tabular}
\end{center}
for some $\beta,\beta' \in [0,4k-1]$.  Now suppose that also $|S_{d,i,2} \cap S_{d',i',2}|=1$.  Since $i \neq i'$, we note that $(4k+2)N-i+1 \neq (4k+2)N-i'+1$.  Also, it cannot be the case that $4kN+4e+i=(4k+2)N-i'+1$, since
\[
4kN+4e+i = (4k+2)N-2d+i \leq (4k+2)N-d < (4k+2)N-d' \leq (4k+2)N-i' < (4k+2)N-i'+1.
\]
Now suppose that $4kN+4e+i=4kN+4e'+i'$.  Then $2d-i=2d'-i'$.  If (a) is true, then $(\beta+2)d=(\beta+2)d'$; since $\beta+2>0$, we have $d=d'$, a contradiction.  On the other hand, if (b) is true, then $(\beta+2)d = (\beta+3)d'$, contradicting Lemma~\ref{squeezy4}.  Thus the only remaining possibility is that $4kN+4e'+i'=(4k+2)N-i+1$, so that $i+i'= 2N-4e'+1=2d'+1$ is odd.  Since $d$ and $d'$ have the same parity, this contradicts (a), so it must be that (b) is true.  It follows that
\[
(\beta+3)d' - \beta d + 1 = 2i \leq 2d.
\]
Thus $(\beta+3)d' \leq (\beta+2)d-1 < (\beta+2)d$, contradicting Lemma~\ref{squeezy4}.

\bigskip

\noindent {\bf Proof of (ii):}  We first note that the largest element in $T_{e,j,1} \cup U_{e,\ell,1}$ is $4kd+(8k-1)e+\ell$, while the smallest element of $T_{e,j,2} \cup U_{e,\ell,2}$ is $4kN+j$.  Since
\[
4kd+(8k-1)e+\ell \leq 4kd+8ke = 4kN < 4kN+j,
\]
it follows that $T_{e,j,1} \cap T_{e',j',2}=\emptyset$, $U_{e,\ell,1} \cap U_{e',\ell',2}=\emptyset$ and $T_{e,j,1} \cap U_{e',\ell',2}=\emptyset$.

Now, if $|T_{e,j,1} \cap T_{e',j',2}| \geq 2$, $|U_{e,\ell,1} \cap U_{e',\ell',1}| \geq 2$ or $|T_{e,j,1} \cap U_{e',j',1}| \geq 2$, then for some $x$, $(x,x+\alpha e) = (x,x+\alpha'e')$, where $\alpha, \alpha' \in [1,8k-1]$.  Thus $\alpha e = \alpha' e'$, contradicting Lemma~\ref{squeezy4}.  If $|T_{e,j,2} \cap T_{e',j',2}| \geq 2$, $|U_{e,\ell,2} \cap U_{e',\ell',2}| \geq 2$ or $|T_{e,\ell,2} \cap U_{e',\ell',2}| \geq 2$, then it follows that $e=e'$, a contradiction.

Thus, if $|T_{e,j} \cap T_{e',j'}| \geq 2$, it must be that $|T_{e,j,1} \cap T_{e',j',1}|=1$ and $|T_{e,j,2} \cap T_{e',j',2}|=1$.  Since $|T_{e,j,1} \cap T_{e',j',1}|=1$, we have that for some
$\alpha, \alpha' \in [0,4k-1]$,
$4kd+\alpha e + j = 4kd' + \alpha'e' +j'$, which implies that $(8k-\alpha)e-j = (8k-\alpha')e'-j'$.  Since $|T_{e,j,2} \cap T_{e',j',2}|=1$, then $4kN+\beta e + j = 4kN+\beta' e' + j'$ where $\beta, \beta' \in \{0,2\}$.
Hence $(8k-\alpha+\beta) e = (8k-\alpha' + \beta') e'$, which contradicts Lemma~\ref{squeezy4} since
$(8k-\alpha+\beta), (8k-\alpha'+\beta') \in [4k+1 , 8k+2]$.
We conclude that $|T_{e,j} \cap T_{e',j'}| \leq 1$.

In a similar way, the assumption that $|U_{e,\ell,1} \cap U_{e,\ell',1}|=1$ and $|U_{e,\ell,2} \cap U_{e,\ell',2}|=1$ leads to a contradiction, as does the assumption that $|T_{e,j,1} \cap U_{e',\ell',1}|=1$ and $|T_{e,j,2} \cap U_{e',\ell',2}|=1$.  We conclude that $|U_{e,\ell} \cap U_{e',\ell'}| \leq 1$ and $|T_{e,j} \cap U_{e',\ell'}| \leq 1$.

Next, suppose that $|U_{e,\ell,1} \cap U_{e',\ell',1}|=1$ and $|U_{e,\ell,2} \cap U_{e',\ell',2}|=1$.  Since $|U_{e,\ell,1} \cap U_{e',\ell',1}|=1$, we have that for some $\alpha, \alpha' \in [4k , 8k-1]$, $4kd+\alpha e + \ell = 4kd' + \alpha' e'+\ell'$, which implies that $(8k-\alpha)e-\ell = (8k-\alpha')e'-\ell$.  Since $|U_{e,\ell,2} \cap U_{e',\ell',2}|=1$, then $4kN+\beta e + \ell = 4kN+\beta'e'+\ell'$ where $\beta, \beta' \in \{1,3\}$.  Hence $(8k-\alpha+\beta)e = (8k-\alpha'+\beta')e'$, which contradicts Lemma~\ref{squeezy4} since $(8k-\alpha+\beta), (8k-\alpha'+\beta') \in [2 , 4k+3]$.

Finally, suppose that $|T_{e,j,1} \cap U_{e',\ell',1}|=1$ and $|T_{e,j,2} \cap U_{e',\ell',2}|=1$.  Since $|T_{e,j,1} \cap U_{e',\ell',1}|=1$, we have that for some $\alpha \in [0,4k-1]$, $\alpha'\in [4k, 8k-1]$, $4kd+\alpha e + j = 4kd' + \alpha'e'+\ell'$, which implies that $(8k-\alpha)e-j = (8k-\alpha')e'-\ell'$.    Since $|T_{e,j,2} \cap U_{e',\ell',2}|=1$, then $4kN+\beta e + j = 4kN+\beta'e'+\ell'$, where
$\beta\in \{0,2\}$ and
$\beta'\in \{1,3\}$.   Hence $(8k-\alpha+\beta)e = (8k-\alpha' + \beta')e'$, which contradicts Lemma~\ref{squeezy4} since $4k+1 \leq 8k-\alpha+\beta\leq 8k+2$ and $2\leq 8k -\alpha' + \beta' \leq 4k+3$.

\bigskip

\noindent {\bf Proof of (iii):}  Note that since
\[
(4k-1)d + i \leq 4kd < 4kN < 4kN+j' < 4kN+e'+\ell',
\]
then $S_{d,i,1} \cap T_{e',j',2} = \emptyset$ and $S_{d,i,1} \cap U_{e',\ell',2} = \emptyset$.  Moreover,
\[
4kd'+(4k-1)e'+j' \leq 4kd' + 4ke' < 4kd' + (8k-1)e' + \ell' \leq 4kd' + 8ke' = 4kN < 4kN+4e+i,
\]
and so $S_{d,i,2} \cap T_{e',j',1} = \emptyset$ and $S_{d,i,2} \cap U_{e',\ell',1} = \emptyset$.

By Lemma~\ref{squeezy4},
\[
(4k-2)d+i \leq (4k-1)d < 4kd'< 4kd'+j'.
\]
It follows that $|S_{d,i,1} \cap T_{e',j',1}| \leq 1$.
Also, since $d<N/3<e'$,
\[
(4k-1)d+i\leq  4kd <4ke' <4kd'+4ke' + \ell',
\]
 and thus $S_{d,i,1} \cap U_{e',\ell',1} = \emptyset$.

Now, using Lemma~\ref{squeezy4}, we also have that
\[
4kN+e'+\ell' \leq 4kN+2e' < 4kN+2e'+j' \leq 4kN+3e' < 4kN+4e < 4kN+4e+i,
\]
and so $S_{d,i,2} \cap T_{e',j',2} = \emptyset$ and $|S_{d,i,2} \cap U_{e',\ell',2}| \leq 1$.  It follows that $|S_{d,i} \cap T_{e',j'}| \leq 1$ and $|S_{d,i} \cap U_{e',\ell'}| \leq 1$.
\end{proof}

\begin{theorem}
Let $k \geq 1$ and $n=(8k+4)N+1$.  There is a set of mutually orthogonal cyclic $(4k+2)$-cycle systems of order $n$ of size at least
$$\frac{N}{96k+30}-1 = \frac{n-1}{(8k+4)(96k+30)}-1.$$  Thus, if $n \equiv 1\pmod{2(4k+2)}$, then
\[
\mu(n,4k+2) \geq \mu'(n,4k+2) \geq \frac{n-1}{(8k+4)(96k+30)} - 1.
\]
\label{Case4k+2}
\end{theorem}

\begin{proof}
The number of integers strictly between $N/3 - N/(48k+15)$ and $N/3$ with the same parity as $N$ is at least $N/(96k+30)-1$.
It thus suffices to show that for distinct integers $d$ and $d'$ with the same parity such that
\[
d,d'\in \left(\frac{N}{3} - \frac{N}{48k+15}, \frac{N}{3}\right),
\]
the decompositions $\mathcal{F}(d)$ and $\mathcal{F}(d')$ are orthogonal.

In turn, it suffices to deal with the exceptional case from Lemma~\ref{Intersections2mod4}.  From Lemma~\ref{specialcycle}, the edges of differences $i$ and $(4k+2)N-i+1$ within $C'(S_{d,i})$ are $\{0,-i\}$ and $\{(4k+2)N-(2k+1)d,-(2k+1)d-i\}$, which are at distance $(4k+2)N - (2k+1)d+i$.  Similarly, the edges of differences $i$ and $(4k+2)N-i+1$ within $C'(S_{d',i})$ are $\{0,-i\}$ and $\{(4k+2)N-(2k+1)d',-(2k+1)d'-i\}$, which are at distance $(4k+2)N-(2k+1)d'+i$.  If the pairs of edges within cycles generated from the starters $C'(S_{d,i})$ and $C'(S_{d',i})$ coincide, then we must have that $(2k+1)d \equiv (2k+1)d'\pmod{n}$.  But $n$ and $2k+1$ are coprime, so $d=d'$.
\end{proof}

\section{Concluding remarks}
\label{conclus}

The main results of this paper have been to establish lower bounds on the number of mutually orthogonal cyclic $\ell$-cycle systems of order $n$.
For upper bounds on the number of systems (not necessarily cyclic in nature) we have the following lemmata.


\begin{lemma} If there exists a set of $\mu$ mutually orthogonal $\ell$-cycle systems of order $n$, then $\mu\leq n-2$.
That is, $\mu(\ell,n)\leq n-2$.
\label{uppa1}
\end{lemma}

\begin{proof}
Consider a vertex $w$ in $K_n$. The vertex $w$ belongs to precisely $(n-1)(n-2)/2$ paths of length $2$ in $K_n$ where $w$ is the center vertex of the path.
Moreover, each such path belongs to at most one $\ell$-cycle from any set of $\mu$ mutually orthogonal $\ell$-cycle systems.
The number of cycles in one $\ell$-cycle system which contain vertex $w$ is equal to $(n-1)/2$.
Thus $\mu(n-1)/2\leq (n-1)(n-2)/2$. The result follows.
\end{proof}

\begin{lemma} Let $\ell\geq 4$.
Then $$\mu(\ell,n)\leq \frac{(n-2)(n-3)}{2(\ell-3)}.$$
\label{uppa2}
\end{lemma}

\begin{proof}
Suppose there exist a set $\{{\mathcal F}_1,{\mathcal F}_2,\dots , {\mathcal F}_{\mu}\}$ of mutually orthogonal $\ell$-cycle systems of $K_n$.
Consider an edge $\{v,w\}$ in $K_n$. Then for each $i\in [1,\mu]$, there is an $\ell$-cycle $C_i\in {\mathcal F}_i$ containing the edge $\{v,w\}$. Let $H$ be the clique of size $n-2$ in $K_n$ not including vertices $v$ and $w$. Then
the intersection of $C_i$ with $H$ is a path $P_i$ with $\ell-3$ edges. Moreover, orthogonality implies that the paths in the set $\{P_i\mid i\in [1,\mu]\}$ are pairwise edge-disjoint.
Thus, $(\ell-3)\mu$ is bounded by the number of edges in $H$; that is, $(\ell-3)\mu\leq (n-2)(n-3)/2$.
\end{proof}

Observe that Lemma \ref{uppa2} improves Lemma \ref{uppa1} only if $\ell >(n+3)/2$.
If $\ell>n/\sqrt{2}$, it is not even possible to find a pair of orthogonal cycle systems, as shown in the following lemma.

\begin{lemma} If $2\ell^2>n(n-1)$ then $\mu(\ell,n)\leq 1.$
\label{uppa3}
\end{lemma}

\begin{proof}
Suppose there exists a pair $\{{\mathcal F}_1,{\mathcal F}_2\}$ of mutually orthogonal $\ell$-cycle systems of $K_n$.
Then ${\mathcal F}_1$ and ${\mathcal F}_2$ each contain $n(n-1)/(2\ell)$ cycles of length $\ell$.
 Let $C$ be a cycle in ${\mathcal F}_1$. By the definition of orthogonality, each edge of $C$ intersects a unique cycle in ${\mathcal F}_2$.
 Thus $\ell \leq n(n-1)/(2\ell)$,  contradiction.
\end{proof}

When the systems are required to be cyclic, Lemma \ref{uppa1} can be slightly improved.

\begin{lemma}\label{CyclicUpperBound} Let $n\geq 4$. If there exists a set of $\mu'$ mutually orthogonal cyclic $\ell$-cycle systems of order $n$, then $\mu'\leq n-3$.
That is, $\mu'(\ell,n)\leq n-3$.
\end{lemma}

\begin{proof}
Since $\mu'(\ell,n)\leq \mu(\ell,n)$, Lemma \ref{uppa1} implies that $\mu'(\ell,n)\leq n-2$. Suppose, for the sake of contradiction that $\mu'(\ell,n)=n-2$.
Thus there exists a set of $n-2$ orthogonal cyclic decompositions of $K_n$ where the vertices are labelled with elements of $\Z_n$. Let $a\in [1,(n-1)/2]$.
Suppose that the path $(-a,0,a)$ of length $2$ does not occur in a cycle from one of these decompositions. Then the total number of paths of length $2$ containing $0$ which appear in one of the cycles is less than $(n-1)(n-2)/2$. However, there are $(n-2)(n-1)/2$ cycles containing vertex $0$, contradicting the condition of orthogonality.

Let $C_a$ be the cycle containing the path $(-a,0,a)$ and let ${\mathcal F}$ be the decomposition of $K_n$ containing $C_a$.
 Since our decomposition is cyclic, there is also a cycle $C'\in {\mathcal F}$ containing $(0,a,2a)$; since $C'$ and $C_a$ share an edge we must have $C'=C_a$.
  Inductively, $C_a=(0,a,2a,\dots )$.
  In particular $C_1=(0,1,2,\dots ,n-1)$ and thus $\ell=n$.
But since $\mu'(\ell,n)=n-2 \geq 2$ and $n>(n-1)/2$, there is a cycle $C''\neq C_a$ in a decomposition ${\mathcal F}'\neq {\mathcal F}$
 containing a repeated difference $a\in [1,(n-1)/2]$.
The cycle $C''$ shares two edges with $C_a$, contradicting the condition of orthogonality.
\end{proof}

It is worth noting that for certain congruencies the upper bound in  Lemma~\ref{CyclicUpperBound} can be made significantly smaller.
For example, if $n \equiv 3\pmod{6}$ then $\mu'(3,n)= 1$, because in this case any cyclic decomposition necessarily contains the cycle $(0,n/3, 2n/3)$.

In the appendix we give computational results for $\mu'(\ell,n)$ when $\ell$ and $n$ are small.
As yet we are unaware of any instances for which the bound of Lemma~\ref{CyclicUpperBound} is tight,
and so we ask if equality ever occurs.

\begin{question}
For which values of $\ell$ and $n$, if any, is $\mu'(\ell,n) = n-3$?
\end{question}


\section*{Acknowledgements}
Authors A.C.\ Burgess and D.A.\ Pike acknowledge research support from NSERC Discovery Grants
RGPIN-2019-04328
and RGPIN-2016-04456, respectively.
Thanks are given to
the Centre for Health Informatics and Analytics of the Faculty of Medicine at Memorial University of Newfoundland
for access to computational resources.

\section*{Appendix}

We computed sets of mutually orthogonal cyclic $\ell$-cycle systems of order $n=2\ell+1$ for small values of $\ell$,
and in so doing we empirically determined or bounded $\mu'(\ell,2\ell+1)$ in these cases.
Recall from Lemma~\ref{CyclicUpperBound} that $\mu'(\ell,2\ell+1) \leq 2\ell-2$.

Note that for any cyclic $\ell$-cycle system of order $2\ell+1$, the cycles of the system comprise a single $\Z_{2\ell+1}$-orbit.
To find sets of mutually orthogonal cyclic $\ell$-cycle systems of order $2\ell+1$, we first determined the orbit
for each possible system and then constructed a graph in which each system is represented as a vertex
and adjacency denotes orthogonality.  Maximum cliques were then sought.  The results for $3 \leq \ell \leq 11$ are summarised in Table~\ref{Table:MuPrime}.
For $9 \leq \ell \leq 11$, we found cliques of order 8 but we do not yet know whether larger cliques exist
(the computational task becomes increasingly challenging as the number of systems grows).

\begin{table}[htbp]
$$
\begin{array}{|c|c|r|c|}
\hline
\ell & n = 2\ell +1 & \mbox{No.~of Cyclic Systems} & \mu'(\ell, 2\ell+1) \\
\hline\hline
3  &  7 &         2 ~~~~~~~ &  2 \\
4  &  9 &         6 ~~~~~~~ &  2 \\
5  & 11 &        24 ~~~~~~~ &  4 \\
6  & 13 &       168 ~~~~~~~ &  5 \\
7  & 15 &      1344 ~~~~~~~ &  8 \\
8  & 17 &     11136 ~~~~~~~ &  8 \\
9  & 19 &    128304 ~~~~~~~ &  \geq 8 \\
10 & 21 &   1504248 ~~~~~~~ &  \geq 8 \\
11 & 23 &  19665040 ~~~~~~~ &  \geq 8 \\
\hline
\end{array}
$$
\caption{Number of mutually orthogonal cyclic $\ell$-cycle systems of order $2\ell+1$}
\label{Table:MuPrime}
\end{table}

We now present examples of the $\Z_{2\ell+1}$-orbits for the sets of mutually orthogonal cyclic $\ell$-cycle systems of order $2\ell+1$
that we found.  Each orbit is represented by the differences that occur on the edges of its cycles,
using notation from Section~\ref{Section:MO4CS}.

\subsection*{$\ell = 3$, $n=7$}

{\small
$[1,2,-3]$,
$[1,-3,2]$
}

\subsection*{$\ell = 4$, $n=9$}

{\small
$[1,-2,-3,4]$,
$[1,4,-3,-2]$
}

\subsection*{$\ell = 5$, $n=11$}

{\small
$[1,-2,4,3,5]$,
$[1,3,-2,5,4]$,
$[1,4,5,-2,3]$,
$[1,5,3,4,-2]$
}

\subsection*{$\ell = 6$, $n=13$}

{\small
$[1,2,3,-4,5,6]$,
$[1,-4,-2,3,-5,-6]$,
$[1,5,3,6,-4,2]$,
$[1,-5,-4,3,-6,-2]$,
$[1,6,3,2,5,-4]$
}

\subsection*{$\ell = 7$, $n=15$}

{\small
$[1,2,6,-4,-7,-3,5]$,
$[1,-2,-3,-5,-4,7,6]$,
$[1,3,4,-2,6,-5,-7]$,
$[1,-3,4,2,-5,-6,7]$,
\\
$[1,5,-3,-7,-4,6,2]$,
$[1,6,7,-4,-5,-3,-2]$,
$[1,7,-6,-5,2,4,-3]$,
$[1,-7,-5,6,-2,4,3]$
}

\subsection*{$\ell = 8$, $n=17$}

{\small
$[1,2,3,4,-6,-7,-5,8]$,
$[1,-2,-3,8,5,-6,4,-7]$,
$[1,3,7,-8,-6,2,-4,5]$,
\\
$[1,-3,-5,6,4,-8,7,-2]$,
$[1,4,5,2,-3,6,-7,-8]$,
$[1,5,-7,6,-8,-3,2,4]$,
\\
$[1,-6,-8,2,5,-4,7,3]$,
$[1,-8,-7,5,-3,4,2,6]$
}

\subsection*{$\ell = 9$, $n=19$}

{\small
$[1,2,3,4,5,-6,-7,9,8]$,
$[1,-2,3,-4,-7,6,8,9,5]$,
$[1,5,8,-3,7,-6,-4,9,2]$,
\\
$[1,-5,-7,-6,-8,-2,-4,3,9]$,
$[1,6,-3,-9,2,-7,4,-5,-8]$,
$[1,-6,8,-7,-3,-5,-2,4,-9]$,
\\
$[1,7,3,6,-8,9,-5,2,4]$,
$[1,9,-7,8,4,3,5,2,-6]$
}

\subsection*{$\ell = 10$, $n=21$}

{\small
$[1,2,3,4,5,-7,6,9,-10,8]$,
$[1,-2,3,-4,5,-6,8,-9,-7,-10]$,
$[1,3,-2,-5,-10,9,-6,-8,4,-7]$,
$[1,-6,-10,7,-3,9,-5,-2,-8,-4]$,
$[1,7,-9,-8,5,6,4,3,2,10]$,
$[1,-8,-5,3,-6,-9,7,-10,2,4]$,
$[1,10,3,6,-7,5,-8,-2,4,9]$,
$[1,-10,3,5,-4,-9,-2,-7,8,-6]$
}

\subsection*{$\ell = 11$, $n=23$}

{\small
$[1,2,3,4,5,6,7,8,9,-10,11]$,
$[1,-2,3,-4,5,-6,7,-11,-10,9,8]$,
\\
$[1,3,2,-4,-5,-11,-6,10,8,-7,9]$,
$[1,-3,10,6,4,7,9,11,8,-2,-5]$,
\\
$[1,-4,8,6,3,5,-9,-2,10,-11,-7]$,
$[1,5,-11,-8,-3,9,-7,2,6,-4,10]$,
\\
$[1,10,-7,-8,3,-5,9,4,6,-11,-2]$,
$[1,-11,5,-7,4,2,-9,-6,8,10,3]$
}

\subsection*{}
Below we present examples of mutually orthogonal cyclic $4$-cycle systems of orders $17$ and $25$;
these are mentioned in Section~\ref{Section:MO4CS}.

\subsection*{$\ell = 4$, $n=17$}

{\small
$\{[1,2,3,-6],[4,-5,-7,8]\}$,
$\{[1,-2,-3,4],[5,8,-7,-6]\}$,
$\{[1,-3,-8,-7],[2,4,5,6]\}$,
\\
$\{[1,4,-7,2],[3,-5,8,-6]\}$,
$\{[1,-4,8,-5],[2,7,-3,-6]\}$,
$\{[1,5,2,-8],[3,-4,-6,7]\}$,
\\
$\{[1,-5,-3,7],[2,-6,8,-4]\}$,
$\{[1,-6,-3,8],[2,-4,7,-5]\}$,
$\{[1,-7,-8,-3],[2,6,5,4]\}$,
\\
$\{[1,-8,-6,-4],[2,5,3,7]\}$
}

\subsection*{$\ell = 4$, $n=25$}

{\small
$\{[1,2,3,-6],[4,-5,12,-11],[7,-8,-9,10]\}$,
$\{[1,-2,-3,4],[5,-6,-7,8],[9,-10,-11,12]\}$,
\\
$\{[1,3,4,-8],[2,5,7,11],[6,-9,-10,-12]\}$,
$\{[1,-3,-4,6],[2,-5,-7,10],[8,-11,-9,12]\}$,
\\
$\{[1,4,2,-7],[3,-5,-9,11],[6,-8,-10,12]\}$,
$\{[1,-4,-2,5],[3,6,7,9],[8,-12,-11,-10]\}$,
\\
$\{[1,5,3,-9],[2,-4,10,-8],[6,12,-7,-11]\}$,
$\{[1,-5,-3,7],[2,6,4,-12],[8,11,-10,-9]\}$,
\\
$\{[1,7,-10,2],[3,-11,-4,12],[5,9,-8,-6]\}$,
$\{[1,-7,-8,-11],[2,12,5,6],[3,10,-4,-9]\}$,
\\
$\{[1,8,4,12],[2,-10,-6,-11],[3,-7,9,-5]\}$,
$\{[1,-8,11,-4],[2,-12,3,7],[5,10,-6,-9]\}$,
\\
$\{[1,-9,-3,11],[2,-6,-4,8],[5,12,-10,-7]\}$,
$\{[1,-10,11,-2],[3,8,-5,-6],[4,-7,12,-9]\}$,
\\
$\{[1,-11,-12,-3],[2,8,10,5],[4,-6,9,-7]\}$,
$\{[1,12,-3,-10],[2,7,-5,-4],[6,11,-9,-8]\}$,
\\
$\{[1,-12,8,3],[2,9,-4,-7],[5,11,-6,-10]\}$
}

\end{document}